\documentclass[11pt,twoside,reqno]{amsart}

\usepackage{orcidlink}
\usepackage{mathtools}
\usepackage{todonotes}
\usepackage{microtype}
\usepackage[noadjust]{cite}
\usepackage[OT1]{fontenc}
\usepackage{type1cm}
\usepackage{amssymb}
\usepackage{dsfont}
\usepackage{enumitem}
\usepackage{comment}
\usepackage{xcolor}
\usepackage{subcaption}
\usepackage{tikz}
\usetikzlibrary{calc}
\usepackage{lmodern}

\usepackage{geometry}
\usepackage{hyperref}
\hypersetup{
  colorlinks=true,
  linkcolor=black,
  anchorcolor=black,
  citecolor=black,
  filecolor=black,      
  menucolor=red,
  runcolor=black,
  urlcolor=black,
}
\usepackage[capitalize,noabbrev,nameinlink]{cleveref}

\makeatletter
\AddToHook{cmd/appendix/before}{\def\cref@section@alias{appendix}}
\makeatother

\theoremstyle{plain}
\newtheorem{theorem}{Theorem}[section]
\newtheorem{maintheorem}{Theorem}

\newtheorem{proposition}[theorem]{Proposition}

\newtheorem{lemma}[theorem]{Lemma}

\theoremstyle{remark}

\theoremstyle{definition}

\newtheorem*{definition*}{Definition}

\newtheorem*{question*}{Question}

\AddToHook{env/proposition/begin}{\crefalias{theorem}{proposition}}
\AddToHook{env/lemma/begin}{\crefalias{theorem}{lemma}}
\AddToHook{env/corollary/begin}{\crefalias{theorem}{corollary}}
\AddToHook{env/claim/begin}{\crefalias{theorem}{claim}}
\AddToHook{env/conjecture/begin}{\crefalias{theorem}{conjecture}}
\AddToHook{env/remark/begin}{\crefalias{theorem}{remakr}}
\AddToHook{env/example/begin}{\crefalias{theorem}{example}}
\AddToHook{env/remark/begin}{\crefalias{theorem}{remark}}

\crefformat{equation}{(#2#1#3)}
\crefformat{figure}{#2Figure #1#3}
\crefformat{enumi}{(#2#1#3)}
\crefformat{theorem}{#2Theorem~#1#3}
\crefformat{lemma}{#2Lemma~#1#3}
\crefformat{proposition}{#2Proposition~#1#3}
\crefformat{corollary}{#2Corollary~#1#3}
\crefformat{claim}{#2Claim~#1#3}
\crefformat{section}{#2Section~#1#3}
\crefformat{subsection}{#2Section~#1#3}
\crefformat{remark}{#2Remark~#1#3}

\numberwithin{equation}{section}

\newcommand{\R}{\mathbb{R}}

\newcommand{\N}{\mathbb{N}}

\newcommand{\mtt}[1]{\mathtt{#1}}

\newcommand{\p}{\mathbf{p}}
\newcommand{\q}{\mathbf{q}}
\newcommand{\x}{\mathbf{x}}
\newcommand{\y}{\mathbf{y}}
\newcommand{\z}{\mathbf{z}}

\renewcommand{\geq}{\geqslant}
\renewcommand{\leq}{\leqslant}

\DeclareMathOperator{\dimloc}{dim_{loc}}

\DeclareMathOperator{\ldimloc}{\underline{dim}_{loc}}

\DeclareMathOperator{\dimh}{dim_H}

\DeclareMathOperator{\dima}{dim_A}

\DeclareMathOperator{\cdima}{\mathcal{C}dim_A}
\DeclareMathOperator{\cdimh}{\mathcal{C}dim_H}

\DeclareMathOperator{\id}{Id}

\renewcommand{\atop}[2]{\genfrac{}{}{0pt}{}{#1}{#2}}

\usepackage{tabularx}
\usepackage{booktabs}

\begin{document}

\title[Conformal dimension of carpets]{Conformal dimension of Bedford--McMullen carpets}

\author{Roope Anttila}
\address[Roope Anttila]
        {University of St Andrews \\ 
         Mathematical Institute\\ 
         St Andrews KY16 9SS\\ 
         Scotland}
\email{ra216@st-andrews.ac.uk}

\subjclass[2020]{} 
\date{\today}
\thanks{This work was supported by the EPSRC, grant no. EP/Z533440/1. The author thanks Sylvester Eriksson-Bique, Jonathan Fraser and Aleksi Pyörälä for many useful discussions on quasisymmetric mappings on fractals.}

\begin{abstract}
We show that the Hausdorff dimension of a Bedford--McMullen carpet with non-uniform fibres can be lowered with a quasisymmetry, provided that the defining iterated function system does not have a full column. We prove this by constructing explicit dimension lowering metrics which are quasisymmetrically equivalent to the Euclidean metric.
\end{abstract}

\maketitle

\section{Introduction}\label{sec:intro}
\emph{Quasisymmetric mappings} are generalisations of bi-Lipshcitz maps, which roughly preserve shapes of sets but may distort their sizes; see \cref{sec:preliminaries} for a precise definition. These mappings originate in geometric function theory, particularly the work of Beurling and Ahlfors \cite{BeurlingAhlfors1956}, but they are also ubiquitous in many seemingly unrelated fields of mathematics. In dynamical systems, one encounters them when studying post-critically finite rational maps, in geometric group theory they arise as boundary extensions of quasi-isometries on hyperbolic groups, and they play an important role when developing analysis on metric spaces with no a priori differentiable structure; we refer to the book of Mackay and Tyson \cite{MackayTyson2010} for an overview of these connections and a comprehensive introduction.

An especially important problem across different fields is characterising spaces up to quasisymmetric equivalence, and in this effort quasisymmetric invariants play a vital role. Some of the most useful quasisymmetric invariants come from fractal geometry, and the notion most relevant to this article is the classical \emph{Hausdorff dimension}, denoted by $\dimh$. Hausdorff dimension itself is invariant under bi-Lipschitz maps, but not under quasisymmetric deformations. In fact, since the identity map from a metric space $(X,d)$ to the snowflaked space $(X,d^{\varepsilon})$ is  a quasisymmetry for any $0<\varepsilon<1$, and since snowflaking the metric changes the Hausdorff dimension by a factor of $1/\varepsilon$, any metric space with positive Hausdorff dimension can be mapped to a space with arbitrarily large Hausdorff dimension by a quasisymmetry \cite[Corollary 1.4.18]{MackayTyson2010}. On the other hand, lowering the Hausdorff dimension by a quasisymmetry is not always possible, for instance, any quasisymmetric image of $\R^d$ always has Hausdorff dimension at least $d$ for purely topological reasons \cite[Theorem 1.4.12]{MackayTyson2010}. This motivates the defintion of the \emph{conformal dimension}, as
\begin{equation}\label{eq:cdimh-def}
    \cdimh X\coloneqq \inf\{\dimh f(X)\colon f\text{ is a quasisymmetry}\}.
\end{equation}
This quasisymmetric invariant has its origins in the work of Pansu \cite{Pansu1989} and while it is theoretically very useful, computing its exact value for a given metric space---or even giving non-trivial upper or lower estimates---is a remarkably difficult problem. For instance, while it is known that the conformal Hausdorff dimension of the Sierpi\'nski carpet is strictly smaller than it's Hausdorff dimension, the precise value remains unknown; see \cite{Kwapisz2020} for the best known upper and lower bounds.

\subsection{Main results}
\begin{figure}
    \includegraphics[width=0.3\linewidth]{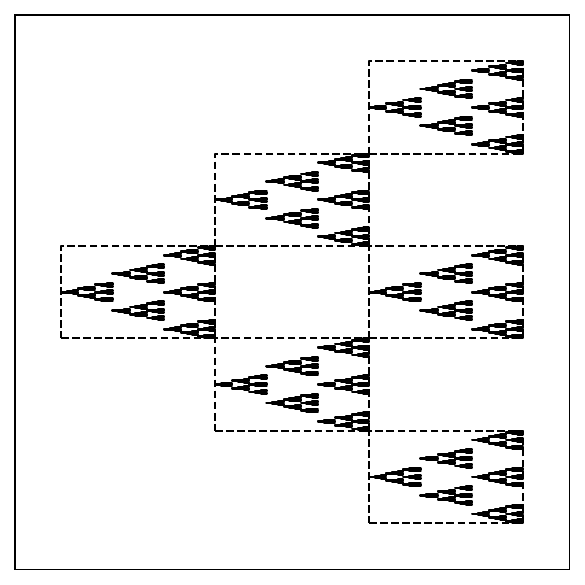}
    \caption{A Bedford--McMullen carpet corresponding to a $3\times 5$ grid. Here $N_1=1$, $N_2=2$ and $N_3=3$, and this example has Hausdorff dimension $\approx1.413$ by \cref{eq:bm-dimh}, and conformal dimension $1$ by \cref{thm:quantitative-ub}.}
    \label{fig:bm}
\end{figure}
Motivated by understanding obstructions in quasisymmetric characterisation of metric spaces, estimating the conformal dimension of fractal spaces has seen growing interest in recent years. Some recent advances include \cite{TysonWu2006,Hakobyan2010,Kwapisz2020,AnttilaEriksson-Bique2026,AnttilaEriksson-BiqueRainio2026} in the deterministic setting, and \cite{RossiSuomala2021,AnttilaEriksson-BiquePyorala2026,MillerTian2026,MillerTian2026a} in the random setting. In this work, we study the conformal dimension of a specific family of deterministic self-affine fractal sets on the plane defined in the following way: Fix integers $1<m<n$ and divide the unit cube into a $m\times n$ grid of rectangles. Choose a subset of these rectangles, and for each chosen rectangle, iterate the proces by dividing said rectangle to an $m\times n$ grid and choosing subrectangles according to the same pattern used initially. The set $X$ which remains after infinitely many steps of this process is called a \emph{Bedford--McMullen carpet}. An example is given in \cref{fig:bm} and a more precise definition as an attractor of a self-affine iterated function system may be found in \cref{sec:preliminaries}. 

The rectangles are naturally divided into $m$ columns, which we enumerate with $i=1,\ldots,m$, and we let $0\leq N_{i}\leq n$ denote the number of chosen rectangles in the $i$th column. Moreover, we let $N_{\min}=\min_{i=1,\ldots,m}N_i$ and $N_{\max}=\max_{i=1,\ldots,m}N_i$ denote the number of rectangles in the smallest and largest column, respectively. In \cite{McMullen1984a}, McMullen proved an analytic formula for the Hausdorff dimension of a Bedford--McMullen carpet in terms of the data which defines the carpet:
\begin{equation}\label{eq:bm-dimh}
    \dimh X=\frac{\log \sum_{i=1}^mN_i^{\frac{\log m}{\log n}}}{\log m}.
\end{equation}
Prior to this work, it was unclear whether the conformal dimension of a Bedford--McMullen carpet is always equal to this value, or if the Hausdorff dimension can be lowered by a quasisymmetry in some cases. The only non-trivial case in which the conformal dimension was known was the case when the Bedford--McMullen carpet has \emph{uniform fibres}, that is when $N_{\max}=N_{\min}\geq 1$. In this case the attractor has a very regular structure, it is in fact Ahlfors regular, and it was shown by Binder, Hakobyan and Li \cite{BinderHakobyanLi2025} using Fuglede modulus that the Hausdorff dimension cannot be lowered by a quasisymmetry.
\begin{theorem}[Binder--Hakobyan--Li]\label{thm:bhl}
    If $X$ is a Bedford--McMullen carpet which satisfies $1\leq N_{\min}=N_{\max}$, then
    \begin{equation*}
        \cdimh X=\dimh X.
    \end{equation*}
\end{theorem}
By choosing a subset of exactly $N_{\min}$ rectangles from each column, any Bedford--McMullen carpet $X$ contains as a subset a Bedford--McMullen carpet $X'$ which has uniform fibres and satisfies
\begin{equation*}
    \dimh X'=1+\frac{\log N_{\min}}{\log n}.
\end{equation*}
Therefore, the result of Binder, Hakobyan and Li gives a lower bound for the conformal dimension of any Bedford--McMullen carpet:
\begin{equation}\label{eq:trivial-lb}
    \cdimh X\geq \cdimh X'\geq 1+\frac{\log N_{\min}}{\log n}.
\end{equation}
However, we note that when $N_{\min}<N_{\max}$, this lower bound is strictly smaller than the Hausdorff dimension in \cref{eq:bm-dimh}. Other than this lower bound, to the best of our knowledge, no non-trivial upper or lower bounds for the conformal dimension were known in the non-uniform fibre case prior to this work. We also note that it was explicitly asked in \cite[Question 8.2]{FraserTyson2026} whether the \cref{thm:bhl} extends to carpets with non-uniform fibres. As our main result, we answer this question in the negative.
\begin{maintheorem}\label{thm:main}
    If $X$ is a Bedford-McMullen carpet which satisfies $1\leq N_{\min}<N_{\max}<n$, then 
    \begin{equation*}
        \cdimh X<\dimh X.
    \end{equation*}
\end{maintheorem}
The condition $N_{\max}<n$ corresponds geometrically to the fact no column is fully included in the construction. This condition is equivalent with $X$ not containing any vertical line segments, and is a technical requirement which we need in the proof. We expect the claim to hold without this assumption as well. We also note that under this condition, the case $N_{\min}=0$ is trivial since the attractor is unifomly disconnected and therefore has $\cdimh X=0$ \cite[Chapter 2.2]{MackayTyson2010}.

\cref{thm:main} follows from a quantitative upper bound for the conformal dimension which only depends on the data defining the carpet, and can be easily computed numerically. Moreover, in the special case when some column has only a single rectangle, our methods show that the conformal dimension is exactly equal to one.

\begin{maintheorem}\label{thm:quantitative-ub}
    Let $X$ be a Bedford-McMullen carpet. If $1=N_{\min}\leq N_{\max}<n$, then
    \begin{equation*}
        \cdimh X=1.
    \end{equation*}
    If $2\leq N_{\min}\leq N_{\max}<n$, then
    \begin{equation*}
        \cdimh X \leq 1+\frac{\log m}{s\log n},
    \end{equation*}
    where $s$ is the unique real number satisfying $\sum_{i=1}^m N_i^{-s}=1$.
\end{maintheorem}

When $2\leq N_{\min}<N_{\max}$, the quantitative upper bound in \cref{thm:quantitative-ub} always satisfies
\begin{equation*}
    1+\frac{\log N_{\min}}{\log n}< 1+\frac{\log m}{s\log n}<\dimh X.
\end{equation*}
Since our upper bound does not match the lower bound in \cref{eq:trivial-lb}, and since our methods do not produce lower bounds for the conformal dimension, the precise value of the conformal dimension unfortunately eludes us when $N_{\min}\geq 2$. 

To conclude this subsection, let us observe how \cref{thm:main} follows form \cref{thm:quantitative-ub}.
\begin{proof}[Proof of \cref{thm:main} assuming \cref{thm:quantitative-ub}]
    In the case $1=N_{\min}<N_{\max}$, it is clear that $\sum_{i=1}^mN_i^{\frac{\log m}{\log n}}>m$, and therefore
    \begin{equation*}
        \cdimh X=1<\frac{\log \sum_{i=1}^mN_i^{\frac{\log m}{\log n}}}{\log m}=\dimh X,
    \end{equation*}
    so \cref{thm:quantitative-ub} implies \cref{thm:main} in this case.

    If $N_{\min}\geq 2$, then we note that since the function $x\mapsto x^{-\frac{\log m}{s\log n}}$ is strictly convex on the positive real line, Jensen's inequality gives
    \begin{equation*}
        \frac{1}{m}\sum_{i=1}^mN_i^{\frac{\log m}{\log n}}=\frac{1}{m}\sum_{i=1}^m(N_i^{-s})^{-\frac{\log m}{s\log n}}\geq\left(\frac{1}{m}\sum_{i=1}^mN_i^{-s}\right)^{-\frac{\log m}{s\log n}}=m^{\frac{\log m}{s\log n}}.
    \end{equation*}
    Moreover, if $N_{\min}<N_{\max}$, then the quantities $N_i^{-s}$ are non-constant, so Jensen's inequality shows that this inequality is strict. By taking logarithms, rearranging and using the Hausdorff dimension formula \cref{eq:bm-dimh}, this gives the claim.
\end{proof}

\subsection{Obstructions for lowering dimension}
Let us briefly discuss why the conclusion of \cref{thm:main} is somewhat surprising. In self-similar settings, such as in the context of the Sierpi\'nski carpet or triangle, constructing dimension lowering quasisymmetries is often possible due to the precise similarity structure of the sets. Informally speaking, if such sets contain large holes at some scale, then they contain large holes of the same shape at all scales due to the self-similar structure, and these holes can be enlarged in a quasisymmetric manner; see \cite{TysonWu2006} or \cite[Chapter 8.3]{MackayTyson2010} for one setting where this can be made precise. For genuinely self-affine sets, such as Bedford--McMullen carpets, there is a crucial difference which makes such constructions more difficult: Even if the set contains large holes at some scale, these holes get flattened by the self-affine distortion at smaller scales, and since quasisymmetries preserve relative shapes, these flattened holes cannot be substantially enlarged in the way described above.

Another way to look at this problem is the fact that, heuristically speaking, the self-affine structure causes the set to look like a product set at small scales, and it is well known that lowering the dimension of product structures with quasisymmetries is difficult, see for example \cite[Proposition 4.1.11]{MackayTyson2010} where it is shown that if $K\subset \R^d$ is the product of a compact set with an interval, then
\begin{equation}\label{eq:conformal-product-dim}
    \cdimh K=\dimh K.
\end{equation}
This has immediate consequences to the quasisymmetric geometry of Bedford--McMullen carpets. Namely, Mackay \cite{Mackay2011} showed that the \emph{Assouad dimension}---a notion of fractal dimension which quantifies the sizes of the thickest parts of a space---cannot be lowered by a quasisymmetry, even in the non-uniform fibre case. Let us denote by $\dima$ the Assouad dimension, see \cite{Fraser2020} for the definition, and by $\cdima$ the \emph{conformal Assouad dimension} defined by replacing Hausdorff dimension with Assouad dimension in \cref{eq:cdimh-def}.
\begin{theorem}[Mackay]\label{thm:mackay}
    Let $X$ be a Bedford--McMullen carpet which satisfies $N_{\min}\geq 1$, then
    \begin{equation*}
        \cdima X=\dima X=1+\frac{\log N_{\max}}{\log n}.
    \end{equation*}
\end{theorem}
Similar results have since been shown to hold for more general self-affine sets, see \cite{AnttilaRutar2024,KaenmakiOjalaRossi2018,BaranyKaenmakiRossi2021}. All of these results are a consequence of \cref{eq:conformal-product-dim} by using the characterisation of the Assouad dimension as the maximal Hausdorff dimension of a \emph{weak tangent}---a limit of a sequence of successive magnifications in the set---after proving that the weak tangent which attains the Assouad dimension is a product of an interval and a compact set.

While the Hausdorff dimension measures the global covering properties of the space instead of local ones, heuristically, this type of product like structure still presists at small scales, which could suggest that the Hausdorff dimension could not be lowered. Our main result shows that this the question for the conformal Hausdorff dimension is more subtle, and this intuition is imprecise.

\subsection{Idea of the proof and structure of the article}
The idea in the proof of \cref{thm:main} is as follows. Observe that the more rectangles a column has, the more it contributes to the Hausdorff dimension of $X$ by \cref{eq:bm-dimh}. Therefore our aim is to construct a metric on $X$, which is quasisymmetric to the Euclidean metric, and shrinks the columns with many rectangles and enlarges the ones with few rectangles. To achieve this, we recall that quasisymmetric metrics on $[0,1]$ are in one-to-one correspondence with doubling measures on $[0,1]$, see \cref{sec:preliminaries} for the precise definitions. We then define a metric on $[0,1]$ corresponding to a suitable doubling measure, which distributes mass to the columns in a manner which is inversely proportional to the number of rectangles in a column. In \cref{sec:metric}, we prove that this metric on $[0,1]$ (in fact \emph{any} metric on $[0,1]$ which is quasisymmetric to the Euclidean metric) extends to a metric on $X$ with uniform multiplicative errors in distances. This result is formalised in \cref{thm:extension-metric}. The construction crucially uses the geometric structure of the carpet $X$, since it is well known that such extensions to the full unit cube do not exist in general. This is also where the assumption $N_{\max}<n$ is used. After constructing the metric, the Hausdorff dimension of $X$ with respect to the new metric may be calculated by essentially the same method which McMullen used in \cite{McMullen1984a} to compute the Hausdorff dimension for the carpet with the usual metric. This is because the diameters of approximate squares in this new metric correspond to the measures of the columns which are easy to compute symbolically. This computation, and the proof of \cref{thm:quantitative-ub}, are concluded in \cref{sec:conformal-dim}.

\section{Preliminaries}\label{sec:preliminaries}
In this section we setup some notation and terminology which will be used throughout the article, and recall some standard concepts from fractal geometry and geometric measure theory.
We start by precisely defining quasisymmetric mappings.

\subsection{Quasisymmetric mappings}
A homeomorphism $f\colon X\to Y$ between metric spaces $(X,d)$ and $(Y,\rho)$ is called an \emph{$\eta$-quasisymmetry} if for all $x,y,z\in X$ with $x\ne z$, we have
\begin{equation*}
    \frac{\rho(f(x),f(y))}{\rho(f(x),f(z))}\leq \eta\left(\frac{d(x,y)}{d(x,z)}\right),
\end{equation*}
where $\eta\colon [0,\infty)\to [0,\infty)$ is an increasing homeomorphism, also called a \emph{distortion function}. A mapping $f\colon X\to Y$ is called a \emph{quasisymmetry} if it is an $\eta$-quasisymmetry for some distortion function $\eta$.

We say that metric spaces $(X,d)$ and $(Y,\rho)$ are \emph{quasisymmetrically equivalent}, if there exists a quasisymmetry $f\colon X\to Y$. It is easy to see that quasisymmetric equivalence is an equivalence relation, see for example \cite[Examples 1.2.10]{MackayTyson2010}, and therefore the conformal dimension defined in \cref{eq:cdimh-def} is invariant under quasisymmetric deformations. If $d$ and $\rho$ are two metrics on the same space $X$, then we say that the metrics $d$ and $\rho$ are \emph{quasisymmetrically equivalent} if the identity function from $(X,d)$ to $(X,\rho)$ is a quasisymmetry, that is, if there exists a distortion function $\eta$, such that
\begin{equation*}
    \frac{\rho(x,y)}{\rho(x,z)}\leq \eta\left(\frac{d(x,y)}{d(x,z)}\right),
\end{equation*}
for all $x,y,z\in X$ with $x\ne z$. If $X\subset \R^d$, then we abuse terminology slightly by saying that a metric $d$ on $X$ is \emph{quasisymmetric}, if it is quasisymmetrically equivalent with the Euclidean metric on $X$.

\subsection{Iterated function systems and self-affine measures}
To define Bedford--McMullen carpets precisely, we recall the concept of iterated function systems. An \emph{iterated function system} is a finite collection $\Phi\coloneqq\{\varphi_i\}_{i\in\Lambda}$ of contractive self-maps of some non-empty compact subset of $\R^d$. Every iterated function system uniquely defines a non-empty and compact \emph{attractor} $X\subset \R^d$, which satisfies
\begin{equation*}
    X=\bigcup_{i\in\Lambda}\varphi_i(X).
\end{equation*}
Similarly, for any fixed probability vector $\p=(p_i)_{i\in\Lambda}$, there exists a unique Borel probability measure $\nu=\nu_{\p}$, which satisfies
\begin{equation*}
    \nu=\sum_{i\in\Lambda}p_i\nu\circ\varphi_{i}^{-1}.
\end{equation*}
If each of the maps $\varphi_i$ is an affine map in $\R^d$, this measure is called the \emph{self-affine measure} corresponding to $\Phi$ and $\p$. In this work, these types of measures are used in two ways, firstly for defining the quasisymmetric metrics on $[0,1]$ which are extended to dimension lowering metrics on $X$, and secondly for computing the Hausdorff dimension of $X$ with respect to these metrics.

It is often convenient to associate an underlying symbolic space to an iterated function system. Given the finite index set $\Lambda$, this associated \emph{symbolic space} is the set $\Lambda^{\N}$. We let $\sigma\colon \Lambda^{\N}\to \Lambda^{\N}$ denote the \emph{left shift} defined by $\sigma(i_1i_2\cdots)=i_2i_3\cdots$. The sets $\Lambda^n$ and $\Lambda^*\coloneqq \bigcup_{n=0}^{\infty}\Lambda^n$ are referred to as \emph{words of length $n$} and \emph{finite words}, respectively. Here $\Lambda^0\coloneqq \{\varnothing\}$ consists of the empty symbol, and we set $\varphi_{\varnothing}=\id$. For $\mtt{i}\in \Lambda^*$, we denote by $|\mtt{i}|$ the \emph{length} of $\mtt{i}$, which is defined as the unique integer which satisfies $\mtt{i}\in\Lambda^{|\mtt{i}|}$. Given $\mtt{i}=i_1\cdots i_{|\mtt{i}|},\mtt{j}=j_1\cdots j_{|\mtt{j}|}\in\Lambda^*$, the \emph{concatenation} of $\mtt{i}$ and $\mtt{j}$ is denoted by $\mtt{ij}\coloneqq i_1\cdots i_{|\mtt{i}|}j_1\cdots j_{|\mtt{j}|}$. For $i\in\Lambda$ and $k\in\N$, we denote by $i^k$, the $k$-fold concatenation of $i$ with itself, that is $i^k=\underbrace{i\cdots i}_{k\text{ times}}$. For $\mtt{i}\in\Lambda^*$, we let
\begin{equation*}
    [\mtt{i}]=\{\mtt{j}\in\Lambda^{\N}\colon \mtt{j}|_{|\mtt{i}|}=\mtt{i}\},
\end{equation*}
denote the \emph{cylinder set} corresponding to $\mtt{i}$. If $\{f_i\}_{i\in\Lambda}$ is any collection of functions indexed by $\Lambda$ and $\mtt{i}=i_1\cdots i_{|\mtt{i}|}\in \Lambda^*$, then we denote by
\begin{equation*}
    f_{\mtt{i}}=f_{i_1}\circ\cdots\circ f_{i_{|\mtt{i}|}},
\end{equation*}
whenever the composition is well defined. Similarly, for any collection $\{a_i\}_{i\in\Lambda}$ of real numbers, we let
\begin{equation*}
    a_{\mtt{i}}=a_{i_1}\cdots a_{i_{|\mtt{i}|}}.
\end{equation*}
Finally, for an infinite word $\mtt{i}=i_1i_2\cdots\in\Lambda^{\N}$ and $n\in\N$, we denote by $\mtt{i}|_n\coloneqq i_1\cdots i_n$ the projection to the first $n$ coordinates.

There is a \emph{natural projection} between the elements of the symbolic space $\Lambda^{\N}$ and points on the attractor $X$. We denote this map by $\Pi\colon \Lambda^{\N}\to X$, and it is defined by setting
\begin{equation*}
    \Pi(\mtt{i})=\lim_{n\to\infty}\varphi_{\mtt{i}|_n}(0),
\end{equation*}
for all $\mtt{i}\in\Lambda^{\N}$. Since each $\varphi_{i}$ is a contraction, it is easy to see that the limit exists for all $\mtt{i}\in\Lambda^{\N}$. Moreover, it is not difficult to see that the map is always surjective, but in general might fail to be injective.

There is also a symbolic characterisation for the self-affine measures $\nu_{\p}$. We denote by $\overline{\nu}=\overline{\nu}_{\q}$ the \emph{Bernoulli measure} corresponding to the probability vector $\q$ on the symbolic space $\Lambda^{\N}$. This is the unique Borel probability measure on $\Lambda^{\N}$, which satisfies
\begin{equation*}
    \overline{\nu}([\mtt{i}])=q_{\mtt{i}},
\end{equation*}
for all $\mtt{i}\in \Lambda^*$. The self-affine measure $\nu$ is then equal to the pushforward of $\overline{\nu}$ under the natural projection $\Pi$, that is $\nu=\overline{\nu}\circ \Pi^{-1}$.

\subsection{Bedford--McMullen carpets}
Bedford--McMullen carpets are now formally defined as a special case of the general construction in the previous section. For $n\in\N$, we denote by $[n]\coloneqq \{1,\ldots,n\}$. Fix integers $n>m>1$ and let
\begin{equation*}
    f_i(x)=\frac{x+i-1}{m}\quad\text{and}\quad g_j(y)=\frac{y+j-1}{n},
\end{equation*}
for each $i\in[m]$ and $j\in[n]$. For any $(i,j)\in[m]\times [n]$, we define $\varphi_{i,j}\colon \R^2\to\R^2$ by setting $\varphi_{i,j}=(f_i,g_j)$. The attractor of the iterated function system $\Phi=\{\varphi_{i,j}\}_{(i,j)\in\Lambda}$, where $\Lambda\subset [m]\times [n]$ is called a \emph{Bedford--McMullen carpet}. For $i\in[m]$, we denote by
\begin{equation*}
    \Lambda_i=\{(i,j)\in\Lambda\colon j\in[n]\},
\end{equation*}
the \emph{column} associated to $i$. Moreover, we let
\begin{equation*}
    N_i=\#\Lambda_i,\quad N_{\min}=\min_{i\in[m]}N_i\quad \text{and}\quad N_{\max}=\max_{i\in[m]}N_i.
\end{equation*}
From the remainder of the article, the notation above is fixed.

We note that for a Bedford--McMullen carpet, the index set is $\Lambda\subset[m]\times [n]$, and thus the corresponding symbolic space $\Lambda^{\N}$ may be naturally identified with a subset of $[m]^{\N}\times [n]^{\N}$ by identifying the word $(i_1,j_1)(i_2,j_2)\cdots\in\Lambda^{\N}$ with $(i_1i_2\cdots,j_1j_2\cdots)\in [m]^{\N}\times [n]^{\N}$. Thus we refer to the elements of $\Lambda^{\N}$ by $(\mtt{i},\mtt{j})$ by implicitly using this natural identification.

Let us next recall the concept of approximate squares, which are very useful when studying the geometry of Bedford--McMullen carpets. For $k\in\N$, we denote by $L(k)$, the unique integer which satisfies
\begin{equation*}
    m^{-(k+1)}\leq n^{-L(k)}<m^{-k}.
\end{equation*}
For $(\mtt{i},\mtt{j})\in\Lambda^{\N}$ and $k\in\N$, we denote by
\begin{equation*}
    \Xi_k(\mtt{i},\mtt{j})=\{(\mtt{i}',\mtt{j}')\in\Lambda^k\colon \mtt{i}=\mtt{i}|_k,\,\mtt{j}'|_{L(k)}=\mtt{j}|_{L(k)}\},
\end{equation*}
the level $k$ \emph{symbolic approximate square} centered at $(\mtt{i},\mtt{j})$, and by
\begin{equation*}
    \mathcal{Q}_k(\mtt{i},\mtt{j})=\bigcup_{(\mtt{i}',\mtt{j}')\in\Xi_k(\mtt{i},\mtt{j})}\varphi_{\mtt{i}',\mtt{j}'}(X),
\end{equation*}
the level $k$ \emph{(geometric) approximate square} centered at $(\mtt{i},\mtt{j})$. The following standard lemma, which follows from the self-affine structure, calculates the measures of approximate squares and will be crucial later on.
\begin{lemma}\label{lemma:approx-square-mass-def}
    Let $\nu$ be a self-affine measure on $X$ corresponding to the probability vector $\q=(q_{i,j})_{(i,j)\in\Lambda}$. Then for any $(\mtt{i},\mtt{j})\in\Lambda^{\N}$ and $k\in\N$, we have
    \begin{equation*}
        \nu(\mathcal{Q}_k(\mtt{i},\mtt{j}))=q_{\mtt{i}|_{L(k)},\mtt{j}|_{L(k)}}Q_{(\sigma^{L(k)}\mtt{i})|_{k-L(k)}},
    \end{equation*}
    where $Q_{i}=\sum_{j\in\Lambda_i}q_{i,j}$.
\end{lemma}

\subsection{Doubling and uniformly perfect measures}
A Borel probability measure $\mu$ supported on $[0,1]$ is said to be \emph{doubling} on $[0,1]$, if there exists a constant $C\geq 1$, such that
\begin{equation*}
    \mu(B(x,2r))\leq C\mu(B(x,r)),
\end{equation*}
for all $x\in[0,1]$ and $r>0$. Note that this definition forces $\mu$ to have full support on $[0,1]$. It is easy to see that the doubling property implies that if $I\subset J$ are (closed) intervals, which satisfy $|J|\leq M|I|$, for some constant $M\geq 1$, then there exists a constant $C\geq 1$ depending only on $M$, such that
\begin{equation*}
    \mu(J)\leq C\mu(I).
\end{equation*}
We also recall that all doubling measures have a finite \emph{Assouad dimension}, which in particular implies that there are $C=C(\mu),\alpha_{\max}=\alpha_{\max}(\mu)>0$ such that for all $0<r\leq R$ and all $x\in[0,1]$, we have
\begin{equation*}
    \frac{\mu(B(x,R))}{\mu(B(x,r))}\leq C\left(\frac{R}{r}\right)^{\alpha_{\max}}.
\end{equation*}

Towards the other extreme, we say that a measure $\mu$ on $[0,1]$ is \emph{uniformly perfect}, if there exist $0<c<1$, such that
\begin{equation*}
    \mu(B(x,r))\leq c\mu(B(x,2r)),
\end{equation*}
for all $x\in X$ and $r>0$. Similarly to the doubling case, there is a notion of dimension which is intimately connected to uniform perfectness. This dimension is called the \emph{lower dimension}, and it implies the existence of constants $C=C(\mu),\alpha_{\min}=\alpha_{\min}(\mu)>0$, such that for all $0<r\leq R$ and all $x\in[0,1]$,
\begin{equation*}
    \frac{\mu(B(x,r))}{\mu(B(x,R))}\leq C\left(\frac{r}{R}\right)^{\alpha_{\min}}.
\end{equation*}

To finish this subsection, we remark that if $\mu$ is a doubling measure on $[0,1]$ with full support, then $\mu$ is automatically uniformly perfect. This is can be seen by using the doubling property to bound the measures of the annuli $B(x,2r)\setminus B(x,r)$. In particular, it follows that if $\mu$ is a doubling measure with full support on $[0,1]$, then there are constants $C,\alpha_{\min},\alpha_{\max}>0$, such that,
\begin{equation}\label{eq:doubling-dim}
    C^{-1}\left(\frac{R}{r}\right)^{\alpha_{\min}}\frac{\mu(B(x,R))}{\mu(B(x,r))}\leq C\left(\frac{R}{r}\right)^{\alpha_{\max}},
\end{equation}
for all $x\in[0,1]$ and $0<r\leq R$. For the remainder of the article, we use the convention that if $\mu$ is a doubling measure with full support on $[0,1]$, then $\alpha_{\min}(\mu)$ and $\alpha_{\max}(\mu)$ denote positive real numbers, such that \cref{eq:doubling-dim} is satisfied. If $\mu$ is clear from the context, we supress it from the notation and denote $\alpha_{\min}(\mu)$ and $\alpha_{\max}(\mu)$ by $\alpha_{\min}$ and $\alpha_{\max}$, respectively.

\subsection{Notation}
For a metric space $(X,d)$ and a subset $A\subset X$, we denote by
\begin{equation*}
    |A|_d\coloneqq \sup_{x,y\in A}d(x,y),
\end{equation*}
the diameter of $A$ in the metric $d$. The open ball with centre $x\in X$ and radius $r>0$ in the metric $d$ is denoted by $B_d(x,r)$.

We denote elements of $\R^2$ by bold lower case letters such as $\x,\y$ or $\z$. We fix the convention that the corresponding non-bold letters with subscripts $1$ and $2$ denote the first and second coordinate, respectively, i.e. $\x=(x_1,x_2)\in \R^2$. The norm $\|\cdot\|$ is fixed throughout the article as the maximum norm in $\R^2$, that is,
\begin{equation*}
    \|\x\|\coloneqq \max\{|x_1|,|x_2|\},
\end{equation*}
for all $\x\in\R^2$. This norm is chosen for convenience and, since all norms on $\R^2$ are (quasisymmetrically) equivalent, replacing it with any other norm on $\R^2$ does not affect our results.

\section{Quasisymmetric metrics on Bedford-McMullen carpets}\label{sec:metric}
In this section, we build a family of quasisymmetric metrics on the Bedford--McMullen carpet $X$. In the following section, we show that certain metrics in this family lower the Hausdorff dimension of the carpet. The construction of these metrics is based on extending quasisymmetric metrics from $[0,1]$ onto $X$, by minimising suitable weights coming from the metric on the line along all chains connecting the points. This construction was inspired by work of Kigami on self-similar spaces \cite[Definition 2.7]{Kigami2014}, but similar ideas have certainly been used long before this, see for example \cite[Chapter 14]{Heinonen2001}.

While extending quasisymmetric metrics from $[0,1]$ to the full unit cube is not necessarily possible---the metric might collapse on vertical lines---the existence of large vertical gaps in Bedfod--McMullen carpets allows one to extend any quasisymmetric metric $d$ on $[0,1]$ to a quasisymmetric metric on $X$, up to uniform multiplicative errors in distances.

For convenience, we formulate the results of this section in terms of doubling measures on $[0,1]$. The equivalence of doubling measures and quasisymmetric metrics on the unit interval is well known \cite[Remark 13.20(b)]{Heinonen2001}: If $\mu$ is a doubling measure on $[0,1]$, then the metric $d$ defined on $[0,1]$ by
\begin{equation*}
    d(x,y)\coloneqq\mu([\min\{x,y\},\max\{x, y\}]),
\end{equation*}
is quasisymmetric. We note that by the doubling property of $\mu$, $d(x,y)\approx \mu(x,|x-y|)$ up to uniform multiplicative errors. A converse result also holds: If $d$ is a quasisymmetric metric on $[0,1]$, then one obtains a doubling measure $\mu$ by extending the measure defined on open intervals $I\subset [0,1]$ by
\begin{equation*}
    \mu(I)=|I|_d,
\end{equation*}
to the Borel $\sigma$-algebra in the usual way. The main result of this section is the following.
\begin{theorem}\label{thm:extension-metric}
    Let $X$ be a Bedford--McMullen carpet which satisfies $1\leq N_{\min}\leq N_{\max}<n$. Then for any doubling measure $\mu$ on $[0,1]$, there exists a metric $d_{\mu}$ on $X$, and a constant $C>0$, such that
    \begin{equation*}
        C^{-1}\mu(x_1,\|\x-\y\|)\leq d_{\mu}(\x,\y)\leq C\mu(x_1,\|\x-\y\|),
    \end{equation*}
    for all $\x,\y\in X$. Moreover, the metric spaces $(X,d_{\mu})$ and $(X,\|\cdot\|)$ are quasisymmetrically equivalent.
\end{theorem}

\subsection{Construction of the metrics}
For the remainder of the section we assume that $1\leq N_{\min}\leq N_{\max}<n$, and we fix a doubling measure $\mu$ on $[0,1]$. We begin the construction of the metric $d_{\mu}$ on $X$ by defining for each pair $\x,\y\in X$,
\begin{equation*}
    \delta(\x,\y)=\mu(B(\min\{x_1, y_1\},\|\x-\y\|)).
\end{equation*}
Notice that $\delta(\x,\y)>0$ if and only if $\x\ne\y$, and that $\delta(\x,\y)=\delta(\y,\x)$. Moreover, since $|x_1-y_1|\leq \|\x-\y\|$, we have $B(x_1,\|\x-\y\|)\subset B(\min\{x_1, y_1\},2\|\x-\y\|)$ and $B(\min\{x_1, y_1\},\|\x-\y\|)\subset B(x_1,2\|\x-\y\|)$, so by using the doubling property of $\mu$, we see that for all $\x,\y\in X$,
\begin{equation}\label{eq:delta-mu-comparison}
    C^{-1}\mu(B(x_1,\|\x-\y\|))\leq \delta(\x,\y)\leq C\mu(B(x_1,\|\x-\y\|)),
\end{equation}
where $C\geq 1$ is the doubling constant of $\mu$.

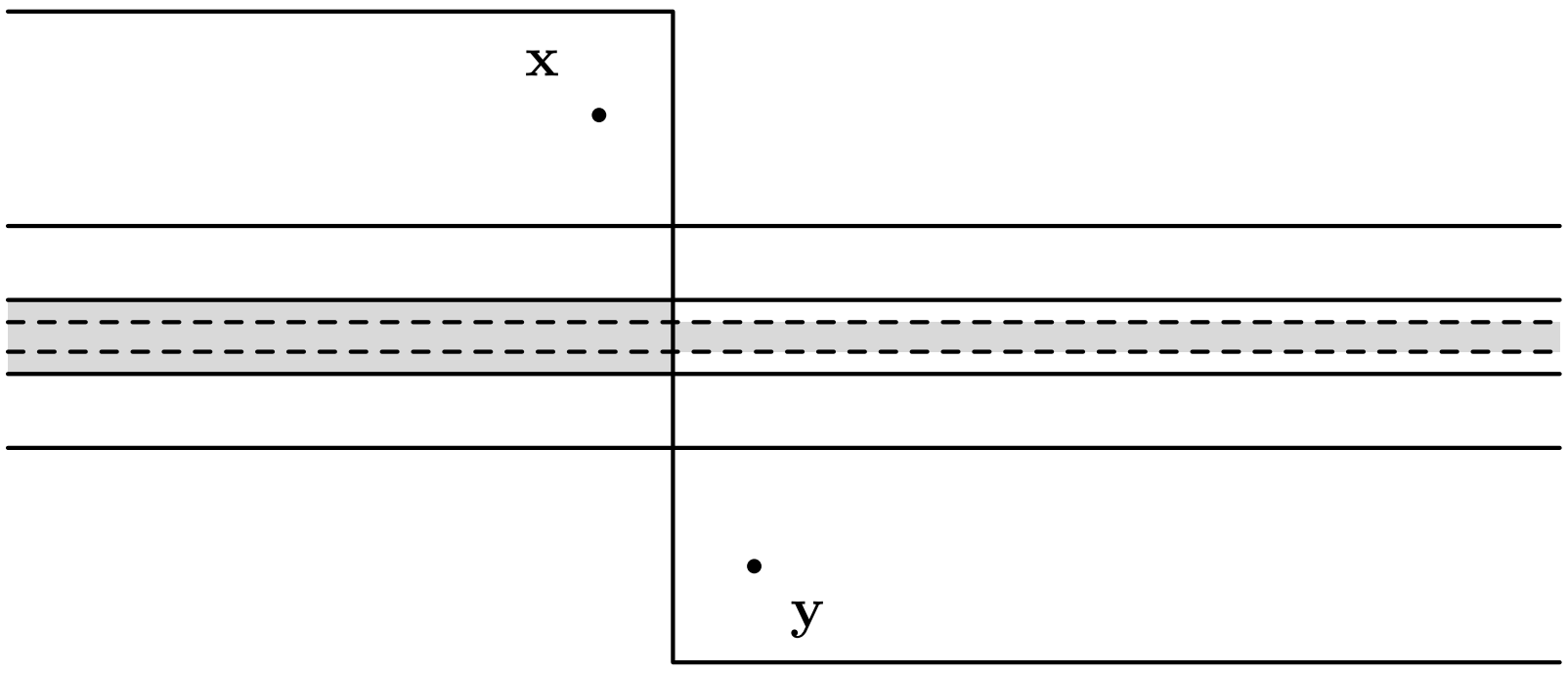
\begin{figure}
    \begin{tikzpicture}[thick,line cap=round,line join=round]

  \def\L{4.5}    
  \def\R{6}      
  \def\H{2.2}    
  \def\h{0.75}   
  \def\a{0.25}   

  \fill[gray!30] (-\L,-\a) rectangle (0,\a);

  \draw (-\L,\a) -- (0,\a);
  \draw (-\L,-\a) -- (0,-\a);

  \def\b{0.10} 

  \fill[gray!30] (0,-\b) rectangle (\R,\b);

  \draw (0,\a) -- (\R,\a);
  \draw (0,-\a) -- (\R,-\a);

  \draw[dashed] (-\L,\b) -- (\R,\b);
  \draw[dashed] (-\L,-\b) -- (\R,-\b);

  \draw (-\L,\h) -- (\R,\h);
  \draw (-\L,-\h) -- (\R,-\h);

  \draw (-\L,\H) -- (0,\H)
                -- (0,-\H) -- (\R,-\H);

  \fill (-0.50,1.5) circle[radius=1.4pt]
    node[above left=5pt] {$\x$};

  \fill (0.55,-1.55) circle[radius=1.4pt]
    node[below right=4pt] {$\y$};

\end{tikzpicture}
    \caption{Illustration of the proof of \cref{lemma:large-gaps} in the worst case when $\x$ is close to the edge of a column. The large half-open rectangles represent cylinders of level $k+1$. One can first find the grey rectangle of level $k+2$ on the left which does not intersect $X$ using the fact that $N_{\max}<n$, and if the corresponding rectangle on the right intersects $X$, one can iterate one step further and use the fact that $N_{\max}<n$ again, to find a rectangle of level $k+3$ which does not intersect $X$.}
    \label{fig:large-gaps}
\end{figure}

We now define the candidate metric $d_{\mu}\colon X\times X\to [0,\infty)$ in the standard way by minimising the $\delta$ value of all chains of points joining $\x$ and $\y$, that is, for $\x,\y\in X$, we set
\begin{equation*}
    d_{\mu}(\x,\y)=\inf\left\{\sum_{j=0}^{M-1}\delta(\x_j,\x_{j+1})\colon \x_0=\x,\,\x_M=\y,\,\x_j\in X,\,\forall j=0,\ldots,M\right\}.
\end{equation*}
Our goal is to prove that the $d_{\mu}$ is indeed a metric on $X$ and satisfies the conclusions of \cref{thm:extension-metric}. We start with a geometric lemma, which formalises the idea that Bedford--McMullen carpets contain large vertical gaps. For $\x\in\R^2$ and $a,b>0$, we denote by
\begin{equation*}
    R_{\x}(a,b)=\{\y\in \R^2\colon |x_1-y_1|< a,\,|x_2-y_2|< b\},
\end{equation*}
the rectangle centred at $\x$ with side lengths $a$ and $b$. The proof of the following lemma is illustrated in \cref{fig:large-gaps}.

\begin{lemma}\label{lemma:large-gaps}
    For any $\x,\y\in X$, there exists $\z\in\R^2$ with $z_1=x_1$, which satisfies
    \begin{equation*}
        R_{\z}(nm^{-5}|x_2-y_2|,n^{-4}|x_2-y_2|)\cap X=\emptyset
    \end{equation*}
\end{lemma}
\begin{proof}
    Fix $\x,\y\in X$ and assume without loss of generality that $x_2<y_2$. Let $k\in\N$ denote the unique integer which satisfies
    \begin{equation}\label{eq:k-def}
        n^{-k}< |x_2-y_2|\leq n^{-(k-1)}.
    \end{equation}
    Let $\mtt{i}_1\in[m]^{k+3}$ be such that $x_1\in f_{\mtt{i}_1}([0,1])$ and let $\mtt{i}_2\in[m]^{k+3}$ denote the coding of the adjacent $m$-adic interval of level $k+3$, chosen such that the distance from $x_1$ to $f_{\mtt{i}_2}([0,1])$ is at most $\tfrac{m^{-(k+3)}}{2}$. We claim that there exists $\mtt{j}\in[n]^{k+3}$, such that $g_{\mtt{j}}([0,1])\subset [x_2,y_2]$, and $(\mtt{i}_1,\mtt{j}),(\mtt{i}_2,\mtt{j})\not \in\Lambda^{k+3}$. Indeed, since $|x_2-y_2|>n^{-k}$, there exists an $n$-adic interval of level $k+1$ between $x_2$ and $y_2$, or in other words, there exists $\mtt{j}_0\in [n]^{k+1}$, such that $g_{\mtt{j}_0}([0,1])\subset [x_2,y_2]$. By using the fact that $N_{\max}<n$, there exists $j_1\in [n]$, such that $(\mtt{i}_1^-,\mtt{j}_0j_1)\not \in\Lambda^{k+2}$ and by using the same fact again if needed, there then exists $j_2\in[n]$, such that  $(\mtt{i}_2,\mtt{j}_0j_1j_2)\not \in\Lambda^{k+3}$. Since $(\mtt{i}_1^-,\mtt{j}_0j_1)\not \in\Lambda^{k+2}$, we also have $(\mtt{i}_1,\mtt{j}_0j_1j_2)\not \in\Lambda^{k+3}$. Setting $\mtt{j}=\mtt{j}_0j_1j_2$, we have $g_{\mtt{j}}([0,1])\subset g_{\mtt{j}_0}([0,1]) \subset [x_2,y_2]$, and $(\mtt{i}_1,\mtt{j}),(\mtt{i}_2,\mtt{j})\not \in\Lambda^{k+3}$.
    
    Since $(\mtt{i}_1,\mtt{j}),(\mtt{i}_2,\mtt{j})\not \in\Lambda^{k+3}$, we have that
    \begin{equation*}
        \varphi_{\mtt{i}_1,\mtt{j}}((0,1)^2)\cup \varphi_{\mtt{i}_2,\mtt{j}}((0,1)^2)\cap X=\emptyset.
    \end{equation*}
    We note that $\varphi_{\mtt{i}_1,\mtt{j}}((0,1)^2)\cup \varphi_{\mtt{i}_2,\mtt{j}}((0,1)^2)$ is a rectangle of width $2m^{-(k+3)}$ and height $n^{-(k+3)}\geq n^{-4}|x_2-y_2|$. Let us denote the $y$-coordinate of the centre line of this rectangle by $y$ and set $\z=(x_1,y)$. Since, by construction, the horizontal distance from $\x$ to either endpoint of $\varphi_{\mtt{i}_1,\mtt{j}}((0,1)^2)\cup \varphi_{\mtt{i}_2,\mtt{j}}((0,1)^2)$ is at least
    \begin{equation*}
        \frac{m^{-(k+3)}}{2}\geq\frac{m^{-3}(n/m)^kn^{-k}}{2}\geq\frac{n}{2m^4} n^{-k}>\frac{n}{m^5}|x_2-y_2|,
    \end{equation*}
    it follows that
    \begin{equation*}
        R_{\z}(nm^{-5}|x_2-y_2|,n^{-4}|x_2-y_2|)\subset \varphi_{\mtt{i}_1,\mtt{j}}((0,1)^2)\cup \varphi_{\mtt{i}_2,\mtt{j}}((0,1)^2),
    \end{equation*}
    and therefore $R_{\z}(nm^{-5}|x_2-y_2|,n^{-4}|x_2-y_2|)\cap X=\emptyset$.
\end{proof}

This following lemma may be used to show both that $d_{\mu}$ is a metric, and that $d_{\mu}$ is quasisymmetric to the Euclidean metric. For an illustration of the proof, see \cref{fig:delta-comparison}

\begin{figure}
    \begin{tikzpicture}[line cap=round, line join=round, thick]

  \def\crosssize{0.045}

  \draw (2.0,-0.65) -- (2.0,2.25);
  \draw (5,-0.65) -- (5,2.25);

\foreach \px/\py in {
  4.15/1.70,
  4.90/1.30,
  5.55/1.20,
  6.20/1.15,
  6.85/0.95,
  7.30/0.30,
  6.70/-0.25,
  6.05/-0.27,
  5.40/-0.40,
  4.80/-0.52
}{
  \fill (\px,\py) circle[radius=0.025];
}

\foreach \px/\py in {
  3.05/1.9,
  2.90/1.35,
  3.20/0.50,
  3.38/0.15,
  3.50/-0.20,
  3.90/-0.58
}{
  \begin{scope}[shift={(\px,\py)}]
    \draw[line width=0.65pt]
      (-\crosssize,-\crosssize)
        -- (\crosssize,\crosssize)
      (-\crosssize,\crosssize)
        -- (\crosssize,-\crosssize);
  \end{scope}
}

  \draw[fill=gray!30] (0,0) rectangle (7,1);

  \fill (3.5,2.25) circle[radius=0.045]
    node[above=5pt] {$\x$};
  \fill (4.20,-0.65) circle[radius=0.045]
    node[below=5pt] {$\y$};

  \draw (-1.05,0) -- (-1.05,1)
    node[midway,left=5pt] {$\approx\|\x-\y\|$};
  \draw (-1.15,0) -- (-0.95,0);
  \draw (-1.15,1) -- (-0.95,1);

  \draw (2.0,-2.05) -- (5,-2.05)
    node[midway,below=5pt] {$\approx\|\x-\y\|$};
  \draw (2.0,-2.15) -- (2.0,-1.95);
  \draw (5,-2.15) -- (5,-1.95);

  \draw (8.1,-0.65) -- (8.1,2.25)
    node[midway,right=6pt] {$\|\x-\y\|$};
  \draw (8.0,-0.65) -- (8.2,-0.65);
  \draw (8.0,2.25) -- (8.2,2.25);

\end{tikzpicture}
    \caption{Illustration of the proof of \cref{lemma:delta-comparison}. By \cref{lemma:large-gaps}, every chain joining $\x$ and $\y$ has to either deviate from $\x$ by a distance comparable to $\|\x-\y\|$ in the horizontal direction (as illustrated by the dotted chain), or cross the vertical gap of height comparable to $\|\x-\y\|$ with a single jump (as illustrated by the crossed chain).}
    \label{fig:delta-comparison}
\end{figure}
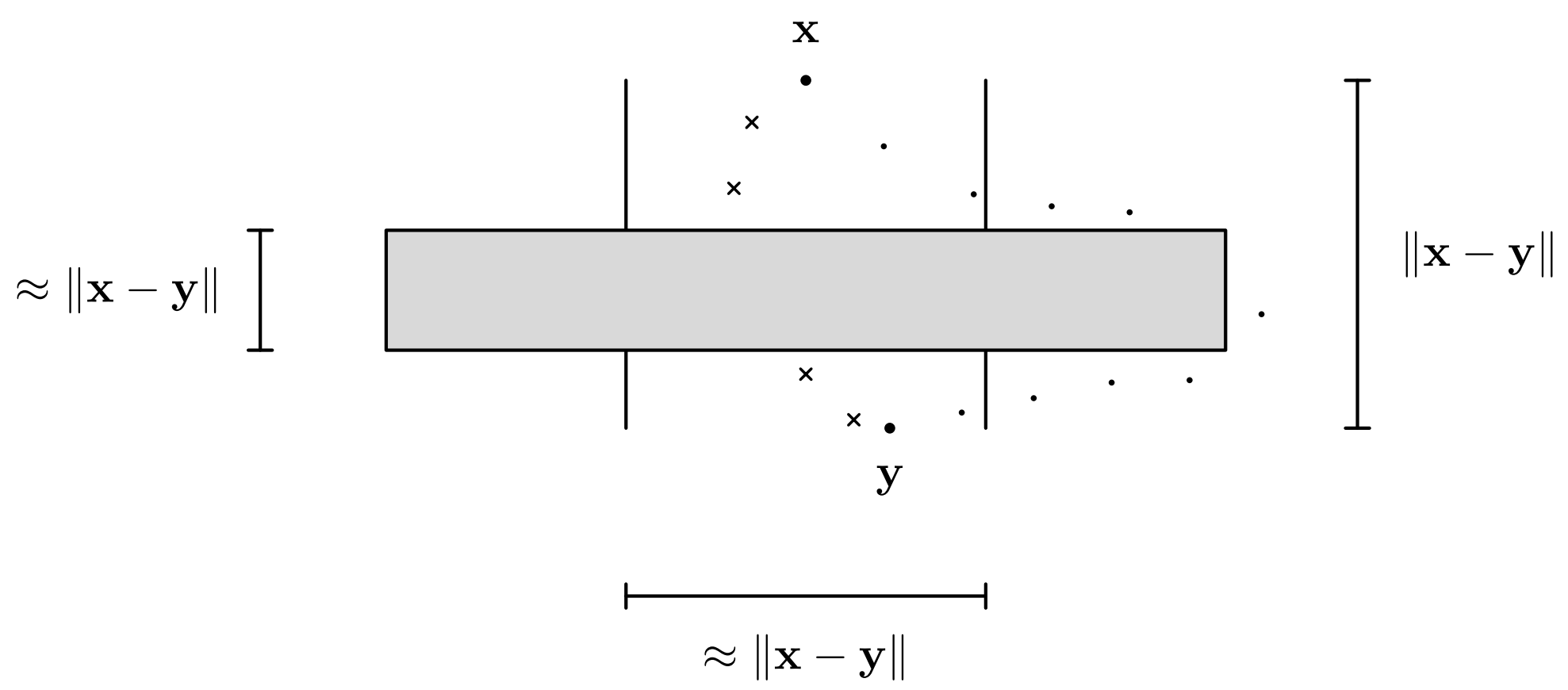
\begin{lemma}\label{lemma:delta-comparison}
    There exists a constant $K\geq 1$, such that for all $\x,\y\in X$, we have
    \begin{equation*}
        K^{-1}\mu(B(x_1,\|\x-\y\|))\leq d_{\mu}(\x,\y)\leq K\mu(B(x_1,\|\x-\y\|)).
    \end{equation*}
\end{lemma}
\begin{proof}
    By considering the chain $\x_0=\x,\x_1=\y$, and using \cref{eq:delta-mu-comparison}, it follows that
    \begin{equation*}
        d_{\mu}(\x,\y)\leq\delta(\x,\y)\leq C\mu(B(x_1,\|\x-\y\|)),
    \end{equation*}
    which proves the upper bound.

    To prove the lower bound, fix $\x,\y\in X$, and assume without loss of generality that $x_1\leq y_1$. Let us also fix an arbitrary chain $\x_0=\x,\x_1,\ldots,\x_M=y$, with $\x_j\in X$ for all $j=0,\ldots,M$. To simplify notation, we denote by $x_j^{\min}=\min\{(x_j)_1, (x_{j+1})_1\}$ and $x_j^{\max}=\max\{(x_j)_1,(x_{j+1})_1\}$. It suffices to show that
    \begin{equation*}
        \sum_{j=0}^{M-1}\delta(\x_j,\x_{j+1})\geq K^{-1}\mu(B(x_1,\|\x-\y\|)),
    \end{equation*}
    for some $K\geq 1$. We split the proof into two cases.
    \newline\newline
    \textbf{Case 1:}
    Assume first that there exists $j_0=0,\ldots,M-1$, such that $|x_1-(x_{j_0})_1|> \tfrac{n^{-4}}{2}\|\x-\y\|$. Since the chain $(\x_j)_{j=0}^{j_0}$ joins $\x$ and $\x_{j_0}$, we have
    \begin{equation*}
        [x_1,(x_{j_0})_1]\subset \bigcup_{j=0}^{j_0-1}B(x_{j}^{\min},|(x_{j+1})_1-(x_{j})_1|)\subset \bigcup_{j=0}^{M-1}B(x_{j}^{\min},\|\x_{j+1}-\x_j\|).
    \end{equation*}
    Indeed, for any $x\in[x_1,(x_{j_0})_1]$, there exists $0\leq j\leq j_0-1$, such that $x\in [x_j^{\min},x_{j}^{\max}]\subset B(x_{j}^{\min},|(x_{j+1})_1-(x_{j})_1|)$ and thus
    \begin{equation*}
        \mu([x_1,(x_{j_0})_1])\leq \sum_{j=0}^{M-1}\delta(\x_j,\x_{j+1}).
    \end{equation*}
    By using the doubling property of $\mu$, we have
    \begin{equation*}
        \mu(B(x_1,\|\x-\y\|))\leq \mu(B(x_1,2n^4|x_1-(x_{j_0})_1|))\leq C_1\mu([x_1,(x_{j_0})_1])\leq C_1\sum_{j=0}^{M-1}\delta(\x_j,\x_{j+1}),
    \end{equation*}
    for some constant $C_1\geq 1$, which only depends on $n$.
    \newline\newline
    \textbf{Case 2:}
    Let us then assume that $|x_1-(x_{j})_1|\leq \tfrac{n^{-4}}{2}\|\x-\y\|$, for all $j=0,\ldots,M$. Since $\x_M=\y$, it follows that
    \begin{equation*}
        |x_1-y_1|\leq \frac{n^{-4}}{2}\|\x-\y\|<\|\x-\y\|.
    \end{equation*}
    This of course implies that $\|\x-\y\|=|x_2-y_2|$. We claim that in this case there exists $j_0=0,\ldots,M-1$, such that $\|\x_{j_0}-\x_{j_0+1}\|\geq n^{-4}\|\x-\y\|$. Indeed, by \cref{lemma:large-gaps}, there is a point $\z\in\R^2$ with $z_1=x_1$, such that
    \begin{equation*}
        R_{\z}(nm^{-5}|x_2-y_2|,n^{-4}|x_2-y_2|)\cap X=\emptyset.
    \end{equation*}
    Since $m<n$, the sides of this rectangle are at horizontal distance at least
    \begin{equation*}
        nm^{-5}|x_2-y_2|>\frac{n^{-4}}{2}|x_2-y_2|,
    \end{equation*}
    from $x_1$ and since the chain $\x_j$ does not deviate by more than $\tfrac{n^{-4}}{2}|x_2-y_2|$ from $x_1$ in the horizontal direction, the chain must cross the vertical gap of height $n^{-4}|x_2-y_2|$ with a single jump at least once, or in other words, there is $j_0=0,\ldots,M-1$, such that
    \begin{equation*}
        \|\x_{j_0}-\x_{j_0+1}\|\geq |(x_{j_0})_2-(x_{j_0+1})_2|\geq n^{-4}|x_2-y_2|=n^{-4}\|\x-\y\|.
    \end{equation*}
    Now, for any $x\in B(x_1,\tfrac{n^{-4}}{2}\|\x-\y\|)$, we have
    \begin{equation*}
        |x_1-(x_{j_0})_1|\leq |x-x_1|+|x_1-(x_{j_0})_1|\leq n^{-4}\|\x-\y\|\leq \|\x_{j_0+1}-\x_{j_0}\|,
    \end{equation*}
    so $B(x_1,\tfrac{n^{-4}}{2}\|\x-\y\|)\subset B((x_{j_0})_1,\|\x_{j_0+1}-\x_{j_0}\|)$. By using \cref{eq:delta-mu-comparison} and the doubling property of $\mu$, we get
    \begin{align*}
        \sum_{j=0}^{M-1}&\delta(\x_j,\x_{j+1})\geq \delta(\x_{j_0},\x_{j_0+1})\geq C^{-1}\mu(B((x_{j_0})_1,\|\x_{j_0}-\x_{j_0+1}\|))\\
        &\geq \mu(B(x_1,\tfrac{n^{-4}}{2}\|\x-\y\|))\geq (CC_2)^{-1}\mu(B(x_1,\|\x-\y\|)),
    \end{align*}
    where $C_2\geq 1$ is some constant which again only depends on $n$. Since $\max\{C,C_1,CC_2\}\leq CC_1C_2$, the claim follows with $K=CC_1C_2$.
\end{proof}
\cref{thm:extension-metric} follows from \cref{lemma:delta-comparison} together with the following proposition.
\begin{proposition}
    The function $d_{\mu}$ is a metric on $X$. Moreover, the metric space $(X,d_{\mu})$ is quasisymmetrically equivalent to $(X,\|\cdot\|)$.
\end{proposition}
\begin{proof}
    It follows from the fact that $\delta(\x,\x)=0$ that $d_{\mu}(\x,\x)=0$, and from $\delta(\x,\y)=\delta(\y,\x)$, that $d_{\mu}(\x,\y)=d_{\mu}(\y,\x)$. Moreover, it is easy to see that $d_{\mu}$ satisfies the triangle inequality, since if $\x,\y,\z\in X$, $\x_0,\ldots,\x_M$ is a chain which joins $\x$ and $\y$ and $\y_0,\ldots,\y_N$ is a chain which joins $\y$ and $\z$, then $\x_0,\ldots,\x_M,\y_1\ldots,\y_N$ is a chain which joins $\x$ and $\z$. Finally, since $\mu(B(x_1,\|\x-\y\|))>0$ if $\x\ne \y$, by applying \cref{lemma:delta-comparison} we get that $d_{\mu}(\x,\y)>0$. Thus $d_{\mu}$ is a metric on $X$.
    
    The fact that $(X,d)$ is quasisymmetrically equivalent with $(X,\|\cdot\|)$ also follows easily from \cref{lemma:delta-comparison}. Let $x,y,z\in X$ with $x\ne z$. Assume first that $\|x-z\|\leq \|x-y\|$. By \cref{lemma:delta-comparison}, and the upper bound in \cref{eq:doubling-dim}
    \begin{equation*}
        \frac{d(x,y)}{d(x,z)}\leq \frac{1}{C^2}\frac{\mu (B(x_1,\|x-y\|))}{\mu(B(x_2,\|x-y\|))}\leq \frac{C'}{C^2}\left(\frac{\|x-y\|}{\|x-z\|}\right)^{\alpha_{\max}},
    \end{equation*}
    where $C'>0$ is the constant of \cref{eq:doubling-dim}. On the other hand, if $\|x-z\|>\|x-y\|$, then by the lower bound in \cref{eq:doubling-dim},
    \begin{equation*}
        \frac{d(x,y)}{d(x,z)}\leq \frac{1}{C^2}\frac{\mu(B(x_1,\|x-y\|))}{\mu(B(x_2,\|x-y\|))}\leq \frac{C'}{C^2}\left(\frac{\|x-y\|}{\|x-z\|}\right)^{\alpha_{\min}}.
    \end{equation*}
    In particular, the identity function is $\eta$-quasisymmetric with $\eta(t)=\tfrac{C'}{C^2} \max\{t^{\alpha_{\min}}, t^{\alpha_{\max}}\}$.
\end{proof}

\section{Upper bound for the conformal dimension}\label{sec:conformal-dim}
In this section, our aim is to prove the quantitative upper bound for the conformal dimension of a Bedford--McMullen carpet stated in \cref{thm:quantitative-ub}. The idea is to construct the dimension lowering metric by the method described in the previous section, by expanding or contracting the columns of the carpet depending on the number of rectangles in each column. Since columns with more rectangles contribute more to the Hausdorff dimension, we want to enlarge the columns with few rectangles and contract the columns with many rectangles. Let us next define the family of measures which we consider.

Let use denote by $\mathcal{P}\coloneqq \{\p\in (0,1)^m\colon \sum_{i=1}^mp_i=1\}$ the simplex of probability vectors in $\R^m$. For a given $\p\in\mathcal{P}$, let $\mu_{\p}$ be the self-affine measure on $[0,1]$ corresponding to the iterated function system $\{f_i\}_{i\in[m]}$ and the probability vector $\p$. It is well known that $\mu_{\p}$ is doubling if and only if $p_1=p_m$, see for example \cite{Yung2007}, so we denote by
\begin{equation*}
    \mathcal{P}_{d}=\left\{\p\in(0,1)^m\colon \sum_{i=1}^mp_i=1\text{ and }p_1=p_m\right\}.
\end{equation*}
For $\p\in\mathcal{P}_d$, let $d_{\p}\coloneqq d_{\mu_{\p}}$ denote the metric on $X$ given by applying \cref{thm:extension-metric} to $\mu_{\p}$. Since $(X,d_{\p})$ is quasisymmetrically equivalent with $(X,\|\cdot\|)$, we trivially have
\begin{equation}\label{eq:doubling-dim-ub}
    \cdimh X\leq \inf_{\p\in \mathcal{P}_d}\dimh (X,d_{\p}).
\end{equation}
The aim of this section is to compute the Hausdorff dimensions of $(X,d_{\p})$, and use the upper bound above to deduce \cref{thm:quantitative-ub}.

\subsection{Hausdorff dimension of $(X,d_{\p})$}
We start by calculating the Hausdorff dimension of $(X,d_{\p})$ for a fixed $\p\in\mathcal{P}_d$. The computation very similar to the classical method used by McMullen to calculate the Hausdorff dimension of $(X,\|\cdot\|)$, namely, we construct a self-affine measure $\nu$ corresponding to a suitable probability vector $\q=(q_{i,j})_{(i,j)\in\Lambda}$, and calculate an upper bound for the pointwise dimension of this measure at all points of the attractor, as well as a corresponding lower bound at $\nu$-almost all points.

Before constructing the measure $\nu$, let us prove a lemma, which compares the $d_{\p}$ diameters of approximate squares to the $\mu_{\p}$ measure of corresponding columns. This lemma also helps explain the choice of the probability vector $\q$.
\begin{lemma}\label{lemma:approx-square-diam}
    There is a constant $C>0$, such that for all $(\mtt{i},\mtt{j})\in\Lambda^{\N}$, we have
    \begin{equation*}
        C^{-1}p_{\mtt{i}|_k}\leq|\mathcal{Q}_k(\mtt{i},\mtt{j})|_{d_{\p}}\leq Cp_{\mtt{i}|_{k}}.
    \end{equation*}
\end{lemma}
\begin{proof}
    Fix $(\mtt{i},\mtt{j})\in\Lambda^{\N}$. By the definition of $\mathcal{Q}_k(\mtt{i},\mtt{j})$, we have
    \begin{equation*}
        \|\x-\y\|\leq m^{-k},
    \end{equation*}
    for all $\x,\y\in\mathcal{Q}_k(\mtt{i},\mtt{j})$. It then follows from \cref{thm:extension-metric}, that
    \begin{equation*}
        d_{\p}(\x,\y)\leq \mu_{\p}(x_1,m^{-k}).
    \end{equation*}
    Noting that $x_1\in f_{\mtt{i}|_k}([0,1])$, and that $|f_{\mtt{i}|_k}([0,1])|=m^{-k}$, we have using the doubling property of $\mu_{\p}$, that
    \begin{equation*}
        d_{\p}(\x,\y)\leq \mu_{\p}(B(x_1,m^{-k}))\leq C_1\mu_{\p}(f_{\mtt{i}|_k}([0,1]))=C_1p_{\mtt{i}|_k},
    \end{equation*}
    for some constant $C_1\geq 1$.

    For the other direction, note that the Euclidean diameter of $\mathcal{Q}_k(\mtt{i},\mtt{j})$ is at least $cm^{-k}$ for some $c>0$, so there are $\x,\y\in X$, such that
    \begin{equation*}
        \|\x-\y\|\geq cm^{-k},
    \end{equation*}
    for some $c>0$. Again, by using \cref{thm:extension-metric} and the doubling property of $\mu_{\p}$, we have
    \begin{equation*}
        d_{\p}(\x,\y)\geq C_2^{-1}\mu_{\p}(f_{\mtt{i}|_k}([0,1]))=C^{-1}p_{\mtt{i}|_k},
    \end{equation*}
    for some constant $C_2\geq 1$. Taking $C=\max\{C_1,C_2\}$ gives the claim.
\end{proof}

Let us now first define the probability vector $\q$, and then explain why this vector gives the correct self-affine measure for calculating the Hausdorff dimension of $(X,d_{\p})$. For any $\p\in\mathcal{P}$, let $t=t(\p)$ denote the unique real number which satisfies
\begin{equation*}
    \sum_{i}N_i^{\frac{\log m}{\log n}}p_i^t=1,
\end{equation*}
and set
\begin{equation}\label{eq:q-def}
    q_{i,j}=N_i^{\frac{\log m}{\log n}-1}p_i^t.
\end{equation}
Since $q_{i,j}$ does not depend on $j$, we have that $Q_i=N_iq_{i,j}=N_i^{\frac{\log m}{\log n}}p_i^t$. The reason for the choice of $\q$ is the following: Since for any $\mtt{i}\in[m]^{\N}$, the sequence $p_{\mtt{i}|_k}$ decreases exponentially, and since the $d_{\p}$-diameter of an approximate square of level $k$ is comparable to $p_{\mtt{i}|_k}$ by \cref{lemma:approx-square-diam}, we may calculate the pointwise dimension of the self-affine measure $\nu$ in the metric $d_{\p}$ along approximate squares. By using \cref{lemma:approx-square-diam} and \cref{lemma:approx-square-mass-def} together with the law of large numbers (or Birkhoff's ergodic theorem) in a standard way, we see that for $\overline{\nu}$-almost every $(\mtt{i},\mtt{j})\in \Lambda^{\N}$,
\begin{align*}
    \lim_{k\to\infty}\frac{\log \nu(\mathcal{Q}_k(\mtt{i},\mtt{j}))}{\log |Q_k(\mtt{i},\mtt{j})|_{d_{\p}}}&=\lim_{k\to \infty} \frac{\log q_{\mtt{i}|_{L(k)},\mtt{j}|_{L(k)}}Q_{(\sigma^{L(k)}\mtt{i})|_{k-L(k)}}}{\log p_{\mtt{i}|_k}}\\
    &=\lim_{k\to \infty} \frac{\log \frac{q_{\mtt{i}|_{L(k)},\mtt{j}|_{L(k)}}}{Q_{\mtt{i}|_{L(k)}}}+\log Q_{\mtt{i}|_{k}}}{\log p_{\mtt{i}|_k}}\\
    &=\frac{\frac{\log n}{\log m}\sum_{i,j}q_{i,j}\log \frac{q_{i,j}}{Q_i}+\sum_{i}Q_i\log Q_i}{\sum_i Q_i\log p_i}.
\end{align*}
The specific vector $\q$ is then chosen in a way which maximises this quantity over all probability vectors.
A simple computation using \cref{eq:q-def} shows that for this specific choice of $\q$,
\begin{equation*}
    \frac{\frac{\log n}{\log m}\sum_{i,j}q_{i,j}\log \frac{q_{i,j}}{Q_i}+\sum_{i}Q_i\log Q_i}{\sum_i Q_i\log p_i}=t,
\end{equation*}
which suggests that $\dimh (X,d_{\p})=t$. The following proposition formalises this non-rigorous argument.
\begin{proposition}\label{prop:p-dim}
    For any $\p\in\mathcal{P}_d$, we have
    \begin{equation*}
        \dimh (X,d_{\p})=t(\p).
    \end{equation*}
\end{proposition}
\begin{proof}
    For any $(\mtt{i},\mtt{j})\in\Lambda$ and $k\in\N$, we have by \cref{lemma:approx-square-mass-def} and the definition of $\q$ and $\mathbf{Q}$, that
    \begin{equation}\label{eq:approx-square-mass}
        \nu(\mathcal{Q}_k(\mtt{i},\mtt{j}))=q_{\mtt{i}|_{L(k)},\mtt{j}|_{L(k)}}Q_{(\sigma^{L(k)}\mtt{i})|_{k-L(k)}}=N_{\mtt{i}|_k}^{\frac{\log m}{\log n}}N_{\mtt{i}|_{L(k)}}^{-1}p_{\mtt{i}|_k}^t.
    \end{equation}
    It follows from \cite[Lemma 4]{McMullen1984a} that for any $\mtt{i}\in\Lambda^{\N}$,
    \begin{equation*}
        \limsup_{k\to\infty} \left(N_{\mtt{i}|_k}^{\frac{\log m}{\log n}}N_{\mtt{i}|_{L(k)}}^{-1}\right)^{\frac{1}{k}}\geq 1.
    \end{equation*}
    Since $p_{\min}^{-\varepsilon}>1$ for any $\varepsilon>0$, there are arbitrarily large $k\in\N$, such that
    \begin{equation*}
        N_{\mtt{i}|_k}^{\frac{\log m}{\log n}}N_{\mtt{i}|_{L(k)}}^{-1}\geq p_{\min}^{-\varepsilon k},
    \end{equation*}
    which together with \cref{eq:approx-square-mass} and the fact that $p_{\mtt{i}|_k}\geq p_{\min}^k$, shows that 
    \begin{equation*}
        \liminf_{k\to\infty}\frac{\log \nu(Q_k(\mtt{i},\mtt{j}))}{\log p_{\mtt{i}|_k}}\leq t-\varepsilon.
    \end{equation*}
    \cref{lemma:approx-square-diam} implies that $\mathcal{Q}_k(\mtt{i},\mtt{j})\subset B_{d_{\p}}(\pi(\mtt{i},\mtt{j}),Cp_{\mtt{i}|_k})$ for some constant $C>0$, which combined with the previous inequality and the fact that the lower local dimension can be calculated along any exponentially decreasing sequence gives
    \begin{equation*}
        \ldimloc(\nu,\pi(\mtt{i},\mtt{j}))=\liminf_{k\to\infty}\frac{\log \nu(B_{d_{\p}}(\pi(\mtt{i},\mtt{j}),Cp_{\mtt{i}|_k}))}{\log Cp_{\mtt{i}|_k}}\leq t-\varepsilon.
    \end{equation*}
    Since $\varepsilon>0$ and $(\mtt{i},\mtt{j})\in\Lambda^{\N}$ were arbitrary, by a standard result in fractal geometry, \cite[Proposition 2.3]{Falconer1997},
    \begin{equation*}
        \dimh(X,d_{\p})\leq t.
    \end{equation*}
    
    Since the upper bound is more important to us, and since the proof of the lower bound is rather standard, we only sketch the details. For a more precise version of the idea in the proof, see for example \cite{McMullen1984a}. It follows from \cref{lemma:approx-square-diam} and \cite[Proposition 3.1]{KaenmakiLiSuomala2016} that for $\overline{\nu}$-almost every $(\mtt{i},\mtt{j})\in \Lambda^{\N}$, one has
    \begin{equation*}
        \ldimloc(\nu,\pi(\mtt{i},\mtt{j}))=\liminf_{k\to \infty}\frac{\log \nu(Q_k(\mtt{i},\mtt{j}))}{\log p_{\mtt{i}|_k}}.
    \end{equation*}
    If we let $n_i(\mtt{i}|_k)=\#\{1\leq j\leq k\colon i_j=i\}$, then by \cref{eq:approx-square-mass}, we have
    \begin{align*}
        \frac{\log \nu(Q_k(\mtt{i},\mtt{j}))}{\log p_{\mtt{i}|_k}}&=t+\frac{\log m}{\log n}\frac{\log N_{\mtt{i}|_k}}{\log p_{\mtt{i}|_k}}-\frac{\log N_{\mtt{i}|_{L(k)}}}{\log p_{\mtt{i}|_k}}\\
        &=t+\frac{\log m}{\log n}\frac{\log\prod_{i=1}^n N_{i}^{n_{i}(\mtt{i}|_k)}}{\log\prod_{i=1}^n p_{i}^{n_{i}(\mtt{i}|_k)}}-\frac{\log\prod_{i=1}^n N_{i}^{n_{i}(\mtt{i}|_{L(k)})}}{\log\prod_{i=1}^n p_{i}^{n_{i}(\mtt{i}|_k)}}\\
        &=t+\frac{\log m}{\log n}\frac{\sum_{i=1}^m n_i(\mtt{i}|_k)\log N_{i}}{\sum_{i=1}^m n_i(\mtt{i}|_k)\log p_{i}}-\frac{\sum_{i=1}^m n_i(\mtt{i}|_{L(k)})\log N_{i}}{\sum_{i=1}^m n_i(\mtt{i}|_k)\log p_{i}},
    \end{align*}
    It follows from the law of large numbers and the definition of $L(k)$, that for $\overline{\nu}$-almost every $(\mtt{i},\mtt{j})$, one has
    \begin{equation*}
        \lim_{k\to\infty}\frac{n_i(\mtt{i}|_k)}{k}=Q_{i}\text{ and }\lim_{k\to\infty}\frac{n_i(\mtt{i}|_{L(k)})}{k}=\frac{\log m}{\log n}Q_i.
    \end{equation*}
    Combining these with the previous inequality shows that for $\overline{\nu}$-almost every $(\mtt{i},\mtt{j})\in\Lambda^{\N}$,
    \begin{align*}
        \dimloc(\nu,\pi(\mtt{i},\mtt{j}))=\lim_{k\to\infty}\frac{\log \nu(Q_k(\mtt{i},\mtt{j}))}{\log p_{\mtt{i}|_k}}&=t,
    \end{align*}
    and therefore by \cite[Proposition 2.3]{Falconer1997}
    \begin{equation*}
        \dimh (X,d_{\p})\geq t.
    \end{equation*}
\end{proof}

\subsection{Quantitative upper bound for conformal dimension}
Next we prove \cref{thm:quantitative-ub}. We start with the following lemma, which computes the infimum of $t(\p)$ over the full simplex $\mathcal{P}$.
\begin{lemma}\label{lemma:t-value}
    If $N_{\min}=1$, then
    \begin{equation*}
        \inf_{\p\in\mathcal{P}}t(\p)=1.
    \end{equation*}
    If $N_{\min}\geq 2$, then
    \begin{equation*}
        \inf_{\p\in\mathcal{P}}t(\p)=1+\frac{\log m}{s\log n},
    \end{equation*}
    where $s$ is the unique real number which satisfies $\sum_{i=1}^mN_i^{-s}=1$.
\end{lemma}
\begin{proof}
    Let us start with the case $N_{\min}>1$. Note that the function $s\mapsto \sum_{i=1}^mN_i^{-s}$ decreases to $0$ as $s$ increases, and that $\sum_{i=1}^m N_i^{-0}=m>1$, so there exists a unique $s>0$, such that
    \begin{equation*}
        \sum_{i=1}^m N_i^{-s}=1.
    \end{equation*}
    Since the function $x\mapsto x^{1+\frac{\log m}{s\log n}}$ is strictly convex, for any $\p\in\mathcal{P}$ we have by Jensen's inequality that
    \begin{align*}
        \sum_{i=1}^{m}N_i^{\frac{\log m}{\log n}}p_i^{1+\frac{\log m}{s\log n}}&=\sum_{i=1}^mN_i^{-s}\left(\frac{p_i}{N_i^{-s}}\right)^{1+\frac{\log m}{s\log n}}\\
        &\geq \left(\sum_{i=1}^mN_i^{-s}\frac{p_i}{N_i^{-s}}\right)^{1+\frac{\log m}{s\log n}}=1.
    \end{align*}
    Since the function $t\mapsto \sum_{i=1}^{m}N_i^{\frac{\log m}{\log n}}p_i^{t}$ is decreasing, this gives $\inf_{\p\in\mathcal{P}}t(\p)\geq 1+\frac{\log m}{s\log n}$.
    By setting $p_i=N_i^{-s}$, the above inequality is an equality, so the infimum is attained.

    Now consider the case $N_{\min}=1$ and let $i_0\in[m]$ be such that $N_{i_0}=1$. For any $\delta>0$ define a probability vector $\p_{\delta}$ defined by setting $p_{i_0}=1-\delta$ and $p_{i}=\frac{\delta}{m-1}$, for all $i\ne i_0$. For any $\varepsilon>0$, we have
    \begin{equation*}
        \sum_{i=1}^mN_i^{\frac{\log m}{\log n}}p_i^{1+\varepsilon}\leq (1-\delta)^{1+\varepsilon}+\frac{n^{\frac{\log m}{\log n}}}{(m-1)^{1+\varepsilon}}\delta^{1+\varepsilon}.
    \end{equation*}
    For a fixed $\varepsilon>0$, a simple Taylor's approximation argument shows that
    \begin{equation*}
        \sum_{i=1}^mN_i^{\frac{\log m}{\log n}}p_i^{1+\varepsilon}\leq 1-(1+\varepsilon)\delta+o(\delta),
    \end{equation*}
    and thus the left hand side is at most $1$ for small enough $\delta$. Hence
    \begin{equation*}
        t(\p_{\delta})\leq 1+\varepsilon,
    \end{equation*}
    and have $\inf_{\p\in\mathcal{P}}t(\p)\leq1$. The lower bound follows trivially by observing that
    \begin{equation*}
        \sum_{i=1}^mN_i^{\frac{\log m}{\log n}}p_i^{1}\geq N_{\min}^{\frac{\log m}{\log n}}=1.
    \end{equation*}

\end{proof}

An observant reader might notice that the numbers in the previous lemma are equal to the quantitative upper bound stated in \cref{thm:quantitative-ub}. To finish the proof of \cref{thm:quantitative-ub}, there is one final slight problem to overcome: \cref{prop:p-dim} and \cref{eq:doubling-dim-ub} only show that
\begin{equation*}
    \cdimh X\leq \inf_{\p\in\mathcal{P}_d}t(\p),
\end{equation*}
but the infimum of $t(\p)$ over $\mathcal{P}$ and $\mathcal{P}_d$ can be different. In fact, one can show that the two infima are equal if and only if $N_1=N_{m}$. However, iterating the system to a deeper level and considering self-affine measures on this iterated system, which is also a Bedford--McMullen carpet, allows us to deduce \cref{thm:quantitative-ub}.
\begin{proof}[Proof of \cref{thm:quantitative-ub}]
    Start by noticing that for any $k\in\N$, the Bedford--McMullen carpet $X$ is the attractor of the diagonal affine iterated function system $\Phi^k=\{\varphi_{\mtt{i},\mtt{j}}\}_{(\mtt{i},\mtt{j})\in\Lambda^k}$. Let us denote by
    \begin{equation*}
        \mathcal{P}^k_d=\left\{\p=(p_{\mtt{i}})_{\mtt{i}\in[m]^k}\colon\sum_{\mtt{i}\in\Lambda^k}p_{\mtt{i}}=1,\,p_{1^k}=p_{m^k}\right\},
    \end{equation*}
    the collection of probability vectors for which the corresponding self-affine measure $\mu^k_{\p}$ on $\{f_{\mtt{i}}\}_{\mtt{i}\in[m]^k}$ is doubling. We claim that it suffices to prove the following: For any $\varepsilon>0$ and $\p\in\mathcal{P}$, there exists $k\in\N$ and $\p_k\in \mathcal{P}^k_d$, such that
    \begin{equation*}
        t_k(\p_k)\leq t(\p)+\varepsilon,
    \end{equation*}
    where $t_k(\p_k)$ is the unique real number which satisfies
    \begin{equation*}
        \sum_{\mtt{i}\in[m]^k}N_{\mtt{i}}^{\frac{\log m}{\log n}}p_{\mtt{i}}^{t_k(\p_k)}=1.
    \end{equation*}
    Indeed, in this case it follows by applying \cref{prop:p-dim} to the iterated system $\Phi^k$ that $\dimh(X,d_{\mu^k_{\p_k}})=t_k(\p_k)$, and therefore
    \begin{equation*}
        \cdimh X\leq \inf_{\mu\text{ doubling on }[0,1]}\dimh (X,d_{\mu})\leq \inf_{\atop{\p\in \mathcal{P}^k_d}{k\in\N}}t_k(\p)\leq\inf_{\p\in \mathcal{P}}t(\p),
    \end{equation*}
    which implies \cref{thm:quantitative-ub} by \cref{lemma:t-value}.

    Let us prove the remaining claim. Let $\varepsilon>0$ and $\p\in\mathcal{P}$. For any $k\in\N$, we denote by $e_k=1-|p_{0^k}-p_{m^k}|$, and consider the probability vector $\p_k=(p'_{\mtt{i}})_{\mtt{i}\in\Lambda^k}$ defined by
    \begin{equation*}
        p'_{0^k}=p'_{m^k}=\frac{\min\{p_{0^k},p_{m^k}\}}{e_k},
    \end{equation*}
    and $p'_{\mtt{i}}=\frac{p_{\mtt{i}}}{e_k}$ for all $\mtt{i}\in[m]^k\setminus\{0^k,m^k\}$. Clearly $\p_k\in\mathcal{P}^k_d$, and
    \begin{align*}
        \sum_{\mtt{i}\in[m]^k}N_{\mtt{i}}^{\frac{\log m}{\log n}}(p'_{\mtt{i}})^{t(\p)+\varepsilon}\leq e_k^{-(t(p)+\varepsilon)}\sum_{\mtt{i}\in [m]^k}N_{\mtt{i}}^{\frac{\log m}{\log n}}p_{\mtt{i}}^{t(\p)+\varepsilon}\leq e_k^{-(t(\p)+\varepsilon)}\left(\sum_{i=1}^mN_i^{\frac{\log m}{\log n}}p_i^{t(\p)+\varepsilon}\right)^k.
    \end{align*}
    It is easy to see that $e_k\to 1$ as $k\to \infty$, and since
    \begin{equation*}
        \sum_{i=1}^mN_i^{\frac{\log m}{\log n}}p_i^{t(\p)+\varepsilon}<1,
    \end{equation*}
    for any fixed $\varepsilon>0$, we have that
    \begin{equation*}
        e_k^{-(t(\p)+\varepsilon)}\left(\sum_{i=1}^mN_i^{\frac{\log m}{\log n}}p_i^{t(\p)+\varepsilon}\right)^k\to 0,
    \end{equation*}
    as $k\to\infty$. In particular, for large enough $k\in\N$,
    \begin{equation*}
         \sum_{\mtt{i}\in[m]^k}N_{\mtt{i}}^{\frac{\log m}{\log n}}(p'_{\mtt{i}})^{t(\p)+\varepsilon}\leq 1,
    \end{equation*}
    which implies that $t_k(\p_k)\leq t(\p)+\varepsilon$, which was the claim.

\end{proof}

\section*{AI disclosure}
ChatGPT-6-Astra was used to solve the standard optimisation problems in \cref{sec:conformal-dim} and to generate tikz code for drawing \cref{fig:large-gaps,fig:delta-comparison} from hand drawn figures of the author. No AI tools were used for generating mathematical ideas or text.

\bibliographystyle{abbrvurl}
\bibliography{bibliography}

\end{document}